\documentclass[a4paper]{article}
\usepackage{amsfonts,amssymb,amscd,amsmath,amsthm,cite,nicefrac,txfonts}
\usepackage{accents}
\usepackage{stmaryrd}
\usepackage[all]{xy} 
\usepackage{epsfig}
\usepackage{color}

\theoremstyle{definition}
\newtheorem{Thm}{Theorem}[section]

\newtheorem{Lem}[Thm]{Lemma}
\newtheorem{Cor}[Thm]{Corollary}
\newtheorem{Def}[Thm]{Definition}

\newcommand{\cplsf}{<\hspace{-0.5ex}\cdot}
\newcommand{\nTs}{\n{T}'}
\newcommand{\myC}{\mathcal{C}}
\DeclareMathOperator{\Leb}{Leb}
\newcommand{\TLeb}{\text{T-Leb}^*}
\newcommand{\fnto}{\ensuremath{\rightarrow}}  

\newcommand{\la}{\ensuremath{\langle}}                 
\newcommand{\ra}{\ensuremath{\rangle}}                 

\newcommand{\U}{\ensuremath{\bigcup}}                     
\newcommand{\pow}{\ensuremath{\mathfrak{P}}}              
\newcommand{\card}[1]{{\ensuremath{\rvert#1\lvert}} }     

\newcommand{\al}[1]{\ensuremath{{\aleph_{#1}}} }          
\newcommand{\om}[1]{\ensuremath{{\omega_{#1}}} }          

\newcommand{\ho}{\ensuremath{^{\omega}}}                  

\newcommand{\ZFC}{\ensuremath{\text{ZFC}} }               
\newcommand{\ZFCx}{\ensuremath{{\text{ZFC}^\ast}} }        
\newcommand{\ZFCxP}{\ensuremath{{\text{ZFC}^\ast_P}} }        
\newcommand{\ZFCxQ}{\ensuremath{{\text{ZFC}^\ast_Q}} }        

\newcommand{\cI}{\ensuremath{{\mathcal I}}} 
\newcommand{\cJ}{\ensuremath{{\mathcal J}}} 

\newcommand{\ex}{{\ensuremath{\exists}} }

\newcommand{\AND}{\ensuremath{\&}}

\newcommand{\THEN}{\ensuremath{\,\rightarrow\,}}
\newcommand{\IFF}{\ensuremath{\,\leftrightarrow\,} }

\newcommand{\DEFEQ}{\ensuremath{\coloneqq}}
\newcommand{\EQDEF}{\ensuremath{\eqqcolon}}

\newcommand{\bQ}{\ensuremath{\mathbb{Q}}}
\newcommand{\bR}{\ensuremath{\mathbb{R}}}
\newcommand{\bB}{\ensuremath{\mathbb{B}}} 
\newcommand{\bC}{\ensuremath{\mathbb{C}}} 

\newcommand{\esm}{\ensuremath{\prec}}  
\newcommand{\vD}{\ensuremath{\vDash}}

\newcommand{\qb}{\ensuremath{\text{``}}}
\newcommand{\qe}{\ensuremath{\text{''}}}

\newcommand{\mD}[2]{\ensuremath{{\underaccent{\tilde}{\Delta}^#1_#2}}}

\newcommand{\forc}{\ensuremath{\Vdash}}

\newcommand{\incomp}{\ensuremath{\perp}}
\newcommand{\comp}{\ensuremath{\parallel}}
\newcommand{\vG}{\ensuremath{\langle G\rangle}}          
\newcommand{\vGP}{\ensuremath{\langle G_P\rangle}}          
\newcommand{\vGR}{\ensuremath{\langle G_R\rangle}}          
\newcommand{\vGQ}{\ensuremath{\langle G_Q\rangle}}          
\newcommand{\n}[1]{\underaccent{\tilde}{#1}}

\newcommand{\hcn}[1]{\n{#1}}

\DeclareMathOperator{\ro}{ro}

\newcommand{\hb}{{\bar{\eta}}} 
\newcommand{\h}{{\hcn{\eta}}}  
\newcommand{\es}{{\eta^\ast}}   
\newcommand{\eot}{{\eta^\otimes}}   

\newcommand{\Ip}{I^+}  
\newcommand{\IpB}{I^+_\text{Borel}}  
\newcommand{\IB}{I_\text{Borel}}
\newcommand{\coI}{\text{co--}I}  

\DeclareMathOperator{\Gen}{Gen}        
\newcommand{\aGen}{\Gen^{\text{abs}}}  

\makeatletter
\def\overbracket#1{\mathop{\vbox{\ialign{##\crcr\noalign{\kern3\p@}
      \downbracketfill\crcr\noalign{\kern3\p@\nointerlineskip}
      $\hfil\displaystyle{#1}\hfil$\crcr}}}\limits}
\def\underbracket#1{\mathop{\vtop{\ialign{##\crcr
      $\hfil\displaystyle{#1}\hfil$\crcr\noalign{\kern3\p@\nointerlineskip}
      \upbracketfill\crcr\noalign{\kern3\p@}}}}\limits}
\def\overparenthesis#1{\mathop{\vbox{\ialign{##\crcr\noalign{\kern3\p@}
      \downparenthfill\crcr\noalign{\kern3\p@\nointerlineskip}
      $\hfil\displaystyle{#1}\hfil$\crcr}}}\limits}
\def\underparenthesis#1{\mathop{\vtop{\ialign{##\crcr
      $\hfil\displaystyle{#1}\hfil$\crcr\noalign{\kern3\p@\nointerlineskip}
      \upparenthfill\crcr\noalign{\kern3\p@}}}}\limits}
\def\downparenthfill{$\m@th\braceld\leaders\vrule\hfill\bracerd$}
\def\upparenthfill{$\m@th\bracelu\leaders\vrule\hfill\braceru$}
\def\upbracketfill{$\m@th\makesm@sh{\llap{\vrule\@height3\p@\@width.7\p@}}%
  \leaders\vrule\@height.7\p@\hfill
  \makesm@sh{\rlap{\vrule\@height3\p@\@width.7\p@}}$}
\def\downbracketfill{$\m@th
  \makesm@sh{\llap{\vrule\@height.7\p@\@depth2.3\p@\@width.7\p@}}%
  \leaders\vrule\@height.7\p@\hfill
  \makesm@sh{\rlap{\vrule\@height.7\p@\@depth2.3\p@\@width.7\p@}}$}
\makeatother

\begin{document}
\begin{center}
{\Large Preserving Preservation}

{\large Jakob Kellner, Saharon Shelah\footnote{Research 
supported by the United States--Israel Binational Science Foundation.
Publication 828.}}
\end{center}


\section*{Abstract}
We present preservation theorems for countable support iteration of nep forcing
notions satisfying ``old reals are not Lebesgue null'' (section
\ref{sec:random}) and ``old reals are not meager'' (section \ref{sec:cohen}).
(Nep is a generalization of Suslin proper.) We also give some results
for general Suslin ccc ideals (the results are summarized in a diagram
on page \pageref{diagram:implications}).

This paper is closely related to  \cite[XVIII, \S3]{Sh:f} and  \cite{Sh:630}. An
introduction to transitive nep forcing and Suslin ccc ideals can be found in
\cite{preservinglittle}.

\tableofcontents

\section{Notation and Basic Results}\label{sec:basic}

In section \ref{sec:nep}, we will use the notion of
nep forcing, as introduced in \cite{Sh:630}.
We will comment on it there.
For the rest of the paper, we only need
some basic facts about proper forcing and Suslin ccc forcing.

In this paper, the notion $N\esm H(\chi)$ always means that $N$ is 
a {\em countable} elementary submodel.

Forcings are written downwards, i.e. $q<p$ means $q$ is a stronger condition
than $p$. Usually, stronger conditions are denoted by symbols lexicographically
bigger than weaker conditions.

\subsection*{Suslin ccc Ideals}

A candidate is a countable transitive model of some suitable fixed
$\ZFCx\subseteq\ZFC$ (see the comments on normal $\ZFCx$
on page \pageref{remark:normal} for more details).
Let $Q$ be a Suslin ccc forcing with an hereditarily countable 
name for a real, $\h$, such that 
$\forc_Q \h\in{\ho\omega}\setminus V$, and such that in all candidates
$\{\llbracket \h(n)=m\rrbracket, n,m\in\omega\}$ generates $\ro(Q)$.
(Such a real is sometimes called ``generic real''. Note that e.g.
the random or Cohen real has this property for the random or Cohen forcing.)

A Suslin ccc ideal is an ideal 
defined from a $Q$ as above in the following way:
$X\in I$ (or: ``$X$ is null'') iff there is a
Borel--set $A$ s.t. $X\subseteq A$ and $\forc_Q\h\notin A$ (where $A$ is
interpreted as a Borel--code evaluated in $V[G]$, not as a set of $V$).\\ $X\in
\Ip$ (or: ``$X$ is positive'')
 means $X\notin I$, and $X$ is $\coI$ (or: ``$X$ has
measure 1'') means ${\ho\omega}\setminus X\in I$.
The set of positive Borel--sets, Borel$\cap \Ip$, is denoted by
$\IpB$, and Borel$\cap I$ is called $\IB$.

For example, if $Q$ is the random algebra then $I$ are the
Lebesgue--Null--sets, if $Q$ is Cohen forcing, then $I$ are the meager
sets.

$\es$  is called $Q$--generic over $M$ ($\es\in\Gen^Q(M)$),
if there is a $G\in V$ $Q$--generic over $M$ s.t. $\h[G]=\es$.
Since $Q$ will usually be fixed, we will just write 
$\Gen(M)$ instead of $\Gen^Q(M)$.

The following can be found e.g. in \cite{preservinglittle}:
\begin{Lem}\label{lem:basic1}
\begin{enumerate}
\item 
$I$ is a $\sigma$--complete ccc ideal containing all singletons, and\\
$\ro(Q)\equiv\text{Borel}/I$ (as a complete Boolean algebra).
\item\label{item:basic1nr3} For $A$ Borel, ``$q\forc \h\in A$'' and ``$A\in I$'' are $\mD12$,
in particular absolute.
\item 
$\Gen(M)= {\ho\omega}\setminus\U\{A^V:\, A \in \IB\cap M\}$.\\
So $\Gen(M)$ is a Borel--set of measure 1.
\end{enumerate}
\end{Lem}

For any Suslin ccc Ideal there is a notion analogous to the
Lebesgue outer measure. Note however that this generalized 
outer measure will be a Borel set, not a real number:

Let $X$ be any set of reals.  A Borel set $B$ is  (a representant of) the outer
measure (o.m.) of $X$, if $B \supset X$, and for all $B'$ s.t. $X\subset
B'\subset B$: $B\setminus B'$ is null.\\
(Note that instead of ``$B\supset X$'' we could use
``$X\setminus B\in I$'' in the definition, that 
makes no difference modulo $I$, since every nullset is contained 
in a Borel nullset).

Clearly, every $X$ has an outer measure (unique modulo $I$);
the outer measure of a Borel--set $A$ is $A$ itself;
the outer measure of a countable union 
is the union of the outer measures; etc

For the Lebesgue ideal, i.e. $Q$=random, the o.m. of
$X$ (according to our definition) is a Borel--set $B$ containing 
$X$ s.t. $\Leb(B)=\Leb^*(X)$, where $\Leb^*$ is the outer
measure according to the usual definition.

For meager, i.e. $Q$=Cohen, the outer measure of a set $X$
is $2^\omega$ minus the union of all basic open sets $C$ s.t.
$C\cap X$ is meager. (This follows from the fact that 
every positive Borel--set contains (mod $I$) a basic open
set and that there are only countable many basic open sets).

\section{Preservation}\label{sec:preservation}

\begin{Def}\label{def:preserving}
Let $X$ be positive with outer measure $B$.
A forcing  $P$:
\begin{description}
\item[preserves positivity of $X$,] if $\forc_P X\in \Ip$
\item[preserves Borel positivity,] if for all positive Borel--sets $A$,
     $P$ preserves the positivity of $A$ (i.e. $\forc_P A^V\in \Ip$).
\item[preserves positivity,] if for all positive $X$, 
     $P$ preserves the positivity of $X$.
\item[preserves outer measure of $X$,] 
	if $\forc_P(B^{V[G]} \text{ is o.m. of }X)$
\item[preserves Borel outer measure,] if for all Borel--sets $A$,
     $P$ preserves the o.m. of $A$ (i.e. $\forc_P A^{V[G]}\text{ is o.m. of }A^V$).
\item[preserves outer measure,] 
     if for all $X$, $P$ preserves the o.m. of $X$.
\end{description}
\end{Def}

With ``preserving positivity (or o.m.) of $V$'' we mean
preservation for $2^\omega$ (or $\bR$ or $\omega^\omega$,
wherever the ideal $I$ lives).

On page \pageref{diagram:implications} there is a diagram of implications
including these notions.

It is clear that preserving o.m. of $X$ 
implies preserving positivity of $X$ (since
being null is absolute for Borel--sets, and the o.m. of $X$ is a
null--set iff $X$ is null).

For all Suslin ccc ideals, preservation of the o.m. of $V$ is
equivalent to preserving Borel o.m.: Let $A$ be a Borel--set in $V$. 
Then in $V[G]$,
the o.m. of $X:=2^\omega\cap V$ is the disjoint union of the o.m. of 
$X\cap A^{V[G]}=A^V$ and the
o.m. of 
$X\setminus A^{V[G]}=(2^\omega\setminus A)^V$.  So if
the o.m. of $A$ decreases, then the outer measure of $V$
decreases.

Another way to characterize Borel o.m. preserving is:
``No positive Borel--set disjoint to $V$ is added''.

If $Q$ is such that in the forcing extension $V'$ of $V$,
$2^\omega\cap V$ has either outer measure $\emptyset$ or
$2^\omega$, then clearly preservation of positivity
of $V$ implies preservation of Borel outer measure.
Note that this is the case for $Q$=random or Cohen.

For positivity, the equivalence of
preservation of $V$ and of all Borel--sets
is not true in general. 
It does hold however if $Q$ satisfies the condition above
(since then preservation of positivity of $V$ implies
even preservation of Borel outer measure).
Another sufficient condition (that is also satisfied
by Lebesgue--null and meager) is that 
$I$
is ``absolutely Borel--homogeneous'':

\begin{Lem}\label{lem:borelhom}
Assume that $P$ preserves positivity of $V$, and that
$Q$ (i.e. $I$) is such that for every $A,B\in\IpB$ there is an 
$A'\subseteq A$, $A'\in \IpB$, and a Borel function $f:A'\fnto B$
such that (in $V[G_P]$) 
for all $X\subseteq B$: $X\in I\THEN X^{-1}\in I$.
Then $P$ preserves positivity of Borel--sets.
\end{Lem}

\begin{proof} (from \cite{Sh:630})
Assume, $P$ makes $B$ null.
Let $J$ be a maximal family of pos. Borel--sets s.t.
for $A',A''\in J$: $A'\cap A''\in I$ and 
there is a $f_{A'}:A'\fnto B$ as in the assumption. Clearly,
$J$ is countable, and its union is $2^\omega$ (mod $I$).
So in $V[G]$, for each $A'\in J$, $A'\cap V$ is null,
since $A'\cap V=f_{A'}^{-1}(B\cap V)$. So $2^\omega\cap V=\bigcup
A'\cap V$ is null.
\end{proof}

Note that the assumption is necessary.  The easy counterexample is the
following: Let $B_0:=\{x\in2^\omega:\, x(0)=0\}$, $B_1:=2^\omega\setminus B_0$.
Let $Q$ add a $r\in 2^\omega$ s.t.  either $r\in B_0$ and $r$ is random or
$r\in B_1$ and $r$ is Cohen.  $\forc_Q r\notin B$ iff $\forc_Q r\notin B\cap
B_0$ and $\forc_Q r\notin B\cap B_1$, i.e. 
iff $B\cap B_0$ is Lebesgue--null and
$B\cap B_1$ is meager. In particular,
$B_0$ and $B_1$ are positive Borel--sets.
Let $P$ be Cohen forcing.  Then $\forc_P  B_0^V\in I$,
but $\forc_P  B_1^V\notin I$.

Borel positivity (or o.m.) preserving
generally (consistently) does not imply positivity preserving,
not even for Cohen or random.
The standard counterexample is the following:
Let $\bB$ be random forcing, and $\bC_\om1$
be the $\al1$ Cohen algebra (which 
adds $\al1$ many Cohens simultaneously).
If $r$ is $\bB$--generic over $V$,  and $(c_i)$
is $\bC_\om1^{V[r]}$--generic over $V[r]$,
then $(c_i)$ is $\bC_\om1^{V}$--generic over $V$ as well.
So $\bB\ast\bC_\om1$ can be factored as $\bC_\om1\ast P$,
where $P$ is ccc. 
In $V[c_i]$, $X=V\cap 2^\omega$ is not meager, but in
$V[c_i][G_P]=V[r][c_i]$ it clearly is.
On the other hand, 
in $V[c_i]$ is non--meager in $V[c_i][G_P]=V[r][c_i]$
(since the $\{c_i\}$ even form a Luzin set).
So $\bC_\om1$ forces that (for $I$=meager)
some ccc forcing $P$ preserves Borel outer measure,
but not positivity. (If Cohen and random are
interchanged, we get an example for $I$=Lebesgue--null).

However, if $P$ is nep (for example if $P$ is Suslin proper), then Borel
positivity preserving {\em does} imply positivity preserving, and Borel outer
measure preserving implies something similar to outer measure preserving, see
Theorem \ref{thm:nep}.

Note that in any case, preservation of positivity (or outer
measure) is trivially preserved by composition of forcings 
(or equivalently: in successor steps of iterations).
How about limit steps?

In this paper, we will restrict ourselves to countable support iterations. Note
that for example for finite support iterations, in all limit steps of countable
cofinalities Cohen reals are added, so preservation of Lebesgue--positivity is
never preserved in finite support iterations.

Preservation of positivity is connected to preservation of 
generic (e.g. random) reals over models:

\begin{Lem}\label{lem:equivpres1}
If $P$ is proper, $X$ positive, then the following are equivalent:
\begin{enumerate}
\item $P$ preserves the positivity of $X$
\item for all 
$N\esm H(\chi)$, $p\in N$
there is an $\eta\in X$ and $q\leq p$ $N$--gen s.t. $q\forc \eta\in\Gen(N[G_P])$
\item for all $p\in P$ there are unbounded (in $2^\omega$) many
$N\esm H(\chi)$ containing $p$ s.t. for some
$\eta\in X$ and $q\leq p$ $N$--gen:  $q\forc \eta\in\Gen(N[G_P])$
\end{enumerate}
\end{Lem}

Here, $A\subseteq \{N\esm H(\chi)\}$ is called unbounded in $2^\omega$,
if for every $x\in 2^\omega$ there is a $N\in A$ s.t. $x\in N$.

\begin{proof}
$1\rightarrow 2$:
Assume $N\esm H(\chi)$, $G$ $V$- and $N$--generic, $p\in G$. In $V[G]$,
$\Gen(N[G])$ is $\coI$, and $X$ is positive, so 
$\Gen(N[G])\cap X$ is nonempty. Now pick a $q$ forcing this. 

$2\rightarrow 3$ is clear. 

$3\rightarrow 1$:
Assume, $p\forc X\subseteq \n{A}\in \IB$. 
Assume, $N\esm H(\chi)$, such that $p,\n{A}\in N$.
If $q\in G$ $V$--generic, then 
$N[G]\vD\qb \n{A}[G]\in I\qe$,
so in $V[G]$ no $\eta\in X$ can be in $\Gen(N[G])$.
\end{proof}

\begin{Lem}\label{lem:equivpres2}
If $P$ is proper, then the following are equivalent:
\begin{enumerate}
\item $P$ preserves positivity
\item For all $N\esm H(\chi)$, there is a
measure--1 Borel--set $A$ s.t.
for all $p\in N$, $\eta \in A$ there is a $q\leq p$ $N$--generic s.t. 
$q\forc \eta\in\Gen^Q(N[G])$. 
\item For all $p$ there are unbounded (in $2^\omega$) many $N\esm H(\chi)$
containing $p$ such that for some 
measure--1--set $A$:
for all $\eta \in A$ there is a $q\leq p$ $N$--generic s.t. 
$q\forc \eta\in\Gen^Q(N[G]) $ 
\end{enumerate}
\end{Lem}

\begin{proof}
$1 \rightarrow 2$:
Since there are only countable many $p$'s in $N$,
it is enough to show that for all
$N$, $p\in N$ there is a set $A$ as in 2.
So pick $N,p$. 
Let $X:=\{\eta:\, \text{ for all }q\leq p\ N\text{--generic, }
        q\forc \eta\notin\Gen(N[G]) \}$. 
We have to show that $X\in I$.
Otherwise (according to Lemma \ref{lem:equivpres1})
there are $q\leq p$ $N$--generic, and $\eta\in X$ s.t.
$q\forc \eta\in\Gen(N[G])$, a contradiction.

$2\rightarrow 3$ is clear.

$3\rightarrow 1$:
Assume $X\in \Ip$, $p\forc X\subseteq \n B\in \IB$. 
Pick $N\esm H(\chi)$ and $A$ s.t $p,\n B\in N$
and $A$ satisfies 3. So for any $\eta\in X\cap A$ there is a 
$q\leq p$ $N$--generic s.t. $q\forc\eta\in \Gen(N[G])$.
But $X\subseteq \n B[G]\in I\cap N[G]$, a contradiction.
\end{proof}

Why are we interested in preservation of generics over models instead
of preservation of positivity? Because in some important cases,
it turns out that preservation of generics is iterable 
(the simplest example is Cohen, see section \ref{sec:cohen}),
while it is not clear how one can show the iterability 
of preservation of positivity directly.

However, to apply the according iteration--theorems, we will
generally need that {\em all} generics are preserved,
not just a measure--1--set as in Lemma \ref{lem:equivpres2}.

It seems that this stronger condition is really necessary, more specific that
the statement ``preservation of Lebesgue--positivity is preserved in countable
support limits of proper forcing iterations'' 
(and the analog statement for meager) is (consistently) false.  A
counterexample seems to be difficult, but we can give a counterexample to the
following (stronger) statement: ``the preservation of positivity of $X$ is
preserved under c.s.i.'s''.  I.e. we can force that there is an
iteration $P_n$  and a
positive set of reals $X$ such that for all $n\in\omega$, $X$ remains positive
after forcing with $P_n$ (it even has o.m. 1), but $P_\omega$ makes $X$ null
(regardless of what limit we take, c.s., f.s., or any other).

The idea is the following (a more precise construction
follows): Let $\bB_\om1$ be the $\al1$ random algebra (which simultaneously
adds $\al1$ many random reals), and $\bC$ the Cohen algebra. Note that $\bC$
makes $V$ null, and $\bB_\om1$ is outer measure preserving and forces 
that the set
of random reals $\{r_\alpha:\, \alpha\in\om1\}$ is an everywhere positive
Sierpinski set.  Let $P$ be the finite support limit $\bB_\om1\ast
\bC\ast\bB_\om1\ast \bC\ast\dots$.  Now factor $P$ the following way: First add
all the randoms, then the first (former) Cohen, the second, the third etc
(these reals are not Cohens anymore, of course).  One would expect that the
first former Cohen will make only the first $\om1$ many randoms null, the
second only the next $\om1$ many, etc. 
So the set of all randoms will become null only
in the limit.

To make that more precise, we will use the following fact:
\begin{Lem}\label{lem:dochkeincounterexample}
Assume, $P_\omega$ is the finite support limit (i.e. union) of $P_0 \cplsf
P_1\cplsf P_2\dots$, and $Q_\omega$ of $Q_0 \cplsf Q_1\cplsf Q_2\dots$. Assume,
$f: P_\omega\fnto Q_\omega$ is s.t. for all $n$\\
(a) $f\restriction P_n: P_n\fnto Q_n$ is complete, and\\
(b) for all $p\in P_{n+1}$, $q\in Q_n$, $r\in P_n$ a reduction of $p$: $f(r)\comp_{P_n}
    q\ \THEN\ f(p) \comp_{P_{n+1}} q$.\\
Then $f: P_\omega\fnto Q_\omega$ is complete.
\end{Lem}
(If $P$ is a subforcing of $P'$, $p\in P'$, then $r\in P$ is called reduction
of $p$ if for all $p'\in P$: $p'\leq r\ \THEN p'\comp p$.  If $P\cplsf P'$, then
there are reductions for all $p\in P'$, and $r$ reduction of $p$ is equivalent
to $r\forc p\in P'/G_P$).
\begin{proof}
It is clear that $f$ preserves $\leq$ and $\incomp$.  Assume
$D\subseteq P_\omega$ is predense, and let $q\in Q$,
i.e. $q\in Q_n$ for some $n$. We have to show that for some $p\in D$, $q\comp f(p)$. 
Let $p'\in P_n$ be a reduction of $q$.  
For some $p\in D$, $p' \comp p$.
$p\in P_m$, wlog $m\geq n$. Set 
$r_m:=p$. 
In $P_{m-1}$ there is a reduction $r_{m-1}$
of $r_m$ s.t. $r_{m-1}\leq p'$ (just take a reduction 
$\tilde r$ of
a $\tilde p\leq p',p$, and let $r_{m-1}\leq \tilde r, p'$).
We continue this construction to get $r_{m-2}$ etc,
until we get $r_n\leq p'\in P_n$ reduction of $r_{n+1}$.
Since  $r_n\leq p'$, and $p'$ is a reduction of $q$,
$f(r_n)\comp q$. Then $f(r_{n+1})\comp q$ by the 
assumption of the lemma ($r_n\in P_n$, $q\in Q_n$,
$r_n$ a reduction of $r_{n+1}$). So continuing this 
up to $m$, we get $f(p)\comp q$.
\end{proof}

Assume in $V$, $S$ is a definition of 
a forcing (i.e. of $p\in S$ and $q\leq_S p$) 
(using arbitrary parameters of $V$).
$S$ is called strongly absolute,
if the following holds: Let $V'$ be a forcing--extension
of $V$. Then $S$ defines a forcing in
$V'$ as well, and ``$p\in S$'', ``$q\leq_S p$'', and ``$\{p_i:\, i\in I\}$
is a max a.c.'' 
are upwards absolute between $V$ and $V'$.

Usually, only ccc forcings will be strongly absolute
(otherwise maximality will not be preserved).
E.g. Mathias forcing (which is a nice, Suslin proper
forcing but not ccc) is not strongly absolute.

On the other hand, every Suslin ccc forcing is clearly strongly absolute. 
Also, (suitable definitions of) 
$\bB_{\kappa}$ or $\bC_\kappa$ (the $\kappa$ Random-
and Cohen--Algebras) are strongly absolute.

If $f_0:\tilde P\fnto \tilde Q$ is complete, and 
$\tilde P$ forces that $\n S$ is strongly absolute,
then clearly $f_0$ can be extended to a complete 
embedding $f_1: P\ast \n S^{V[G_{\tilde P}]}\fnto Q\ast \n S^{V[G_{\tilde Q}]}$:
Just define $f_1(p,\n \tau):=(f_0(p),f_0^*\n \tau)$,
(where $f_0^*\n \tau$ is a $Q$--name s.t.
$f_0^*\n \tau[G_Q]_Q=\n \tau[f_0^{-1} G_Q]_P$).

Note that $f_1$ is not only complete, but satisfies the
second condition of Lemma \ref{lem:dochkeincounterexample}
as well:
if $r$ is a reduction of $(p,\n \tau)$
(wlog $r=p$), and if $f_0(r)$ is compatible with some $q\in Q$
(wlog $q\leq f_0(r)$), then $f_1(p,\n \tau)$ is compatible with
$q$ by absoluteness.

Therefore we can iterate the extension of $f_0$ and get the 
following:

\begin{Lem}\label{lem:nicelimit}
Let $f:\tilde P\fnto \tilde Q$ be complete, and 
$(R_n,{\n S}_n)_{n\in\omega}$ be (the definition for a) finite support iteration,
and $P\ast R_n$ forces that ${\n S}_n$ is strongly absolute. 
Then $f$ can be extended to a complete embedding of
$\tilde P\ast R_\omega^{V[G_{\tilde P}]} \fnto \tilde Q\ast R_\omega^{V[G_{\tilde Q}]}$.
\end{Lem}

Now we can finally construct the counterexample:
Define $P_n$ to be the finite support limit (at $\omega$) of: first n copies of
$\bB_\om1\ast \bC$, then $\omega$ copies of $\bB_\om1$.  To be able to refer to
the random reals added by $P_n$, we denote the $i$--th copy of $\bB_\om1$ with
$\bB^i_\om1$, and the random reals added by this copy with $r^i_\alpha$
($\alpha\in\om1$). So
$P_n:=\bB^0_\om1\ast\bC\ast \dots\ast \bB^{n-1}_\om1\ast\bC\ast \bB^n_\om1\ast\bB^{n+1}_\om1\ast\dots$.

We claim that there is a complete embedding $f$\\
\begin{tabular}{@{}lllr@{}l@{}}
from & 
  $P_n$ & $=$ &
    $\overbracket{\bB^0_\om1\ast\bC\ast \dots\ast \bB^{n-1}_\om1\ast\bC\ast\bB^n_\om1\ast\phantom{\bC\ast}}\,$&
      $\bB^{n+1}_\om1\ast \bB^{n+2}_\om1\ast\dots$\\ 
to & 
  $P_{n+1}$ & $=$ & 
    $\underbracket{\bB^0_\om1\ast\bC\ast \dots\ast \bB^{n-1}_\om1\ast\bC\ast\bB^n_\om1\ast\bC\ast}\,$&
      $\bB^{n+1}_\om1\ast \bB^{n+2}_\om1\ast\dots$\\
\end{tabular}\\
Lets call the blocks marked above $\tilde P$ and $\tilde Q$, resp.
It is trivial that we find a complete embedding $f: \tilde P\fnto \tilde Q$. 
So by the last lemma, we can extend it to a complete embedding $P_n\fnto P_{n+1}$.
It is also clear that $f$
leads to the same evaluation of the random reals, i.e. it
has the following property: If $G_{n+1}$ is $P_{n+1}$--generic,
and $G_n:=f^{-1} G_{n+1}$ is the corresponding $P_n$--generic
filter, then $r^m_\alpha[G_n]_{P_n}=r^m_\alpha[G_{n+1}]_{P_{n+1}}$
for all $l\in\omega,\alpha\in\om1$.

$P_n$ forces that $\{r^l_\alpha:\, l<n,\alpha\in\om1\}$ is a null--set and that
$\{r^n_\alpha:\, \alpha\in\om1\}$ is not null (it even has outer measure 1).
So in $V[G_0]$ (after forcing with $P_0$), we have a positive set
$X:=\{r^l_\alpha:\, l\in\omega,\alpha\in\om1\}$, and ccc forcings $P_1\cplsf
P_2\cplsf\dots$ such that $X$ has outer measure 1 after forcing with each
$P_n$, but any forcing that adds  generics for all the $P_n$ makes $X$ null
(since $X$ is the countable union of the $\{r^l_\alpha:\, \alpha\in\om1\}$).

Notes:\\
So we (consistently) get a counterexample for the following statement:
The $\omega$--limit of ccc forcings preserving the outer 
measure of $X$ preserves the positivity of $X$.\\
The dual example shows that the preservation of Baire--positivity of
a specific set is consistently not preserved at $\omega$--limits 
(of any iteration).

\section{True Preservation}\label{sec:true}

Preservation of all generics (not just a measure--1--set of them)
is closely related to preserving ``true positivity'', a notion
using the stationary ideal
on $\mathfrak{P}_{\aleph_1}$.

First we recall some basic facts:

\begin{Lem}\label{lem:basicstat}
Let $\cI$ and $\cJ_1\subseteq \cJ_2$ be arbitrary. 
\begin{enumerate}
\item The club--filter on $[\cI]^\al0$ is closed under countable 
	intersections.
\item If $C\subseteq [\cJ_1]^\al0$ is club, then $C^{\cJ_2}:=\{B\in [\cJ_2]^\al0:\, 
  B\cap \cJ_1\in C\}$ is club.
\item If $C\subseteq [\cJ_2]^\al0$ is club, then $C^{\cJ_1}:=\{B\cap \cJ_1:\, 
  B\in C\}$ is club.
\item A forcing $P$ is proper iff for arbitrary $\cI$, and $S\subseteq [\cI]^\al0$
  stationary: $\forc_P S$ is stationary.
\item If $C\subseteq [\cI]^\al0$ is club, and $P$ is proper,
  then in $V[G]$ there is a $C'\subseteq [\cI]^\al0$ club
  s.t. $C'\cap V=C$. 
\end{enumerate}
\end{Lem}
(Note that if $C$ is club in $V$, then generally $C$ will not be club
any more in $V[G]$. To prove the last item,
use the usual basis--theorem for club--sets).

Assume $\cI$ is an arbitrary index-set,
$S\subseteq [\cI]^\al0$ stationary, 
 $\hb=(\eta_s: s\in \cI)$ a sequence of reals.
Pick any $\cJ\supset \cI\cup 2^\omega$.
  For $C\subseteq [\cJ]^\al0$, define\\
$S(C):=\{s\in S:\, \exists N\in C:\, N\cap \cI=s\,\&\, \eta_s\in\Gen^Q(N)\}$\\ 
$\hb(C):=\{\eta_s:\, s\in S(C)\}$.\\
$\eta_s\in\Gen^Q(N)$ means that $\eta_s\notin B$ for all Borel--null--sets $B$
coded by a real in $N$. If $N\esm H(\chi)$ for some regular $\chi$ (and we will
only be interested in this case), then
$\eta_s\in\Gen^Q(N)$ is equivalent to the following:
there is a $G\in V$ $Q$--generic over $N$ s.t. $\eta_s=\h[G]$
(to see this, just apply Lemma \ref{lem:basic1}(\ref{item:basic1nr3}) 
to the transitive collapse of $N$).
  
\begin{Def}
\begin{enumerate}
\item $\hb$ is truly positive,  
  if for all $C\subseteq [\cJ]^\al0$ club, $\hb(C)\in \Ip$. 
\item $B$ is the true outer measure of $\hb$,
  if it is the smallest Borel--set containing any of the $\hb(C)$, i.e.
  if the following holds: $B$ is Borel, for some $C\subseteq [\cJ]^\al0$ club
  $\hb(C)\subseteq B$, and for any
  $C'\subseteq [\cJ]^\al0$  club, $B'\supseteq \hb(C')$: $B\setminus B'\in I$.
\end{enumerate}
\end{Def}

\begin{Lem}\label{lem:truebasic}
\begin{enumerate}
\item the above notions do not depend on $\cJ$ (provided that $\cJ\supset \cI
  \cup 2^\omega$).
\item The true outer measure always exist.
\item  If $\cJ=H(\chi)$, then TFAE:\\
(a) $\hb$ is truly positive\\
(b) for all $C\subseteq [\cJ]^\al0$ club, $\hb(C)\neq\emptyset$\\
(c) for all $x$, there is an $N\esm H(\chi)$ containing $x,\cI,S,\hb$ s.t.
$N\cap I=s\in S$ and $\eta_s\in\Gen^Q(N)$.
\end{enumerate}
\end{Lem}

\begin{proof}
\begin{enumerate}
\item Assume, $\cI\cup 2^\omega\subseteq \cJ_1\subseteq \cJ_2$. 
  Assume, $C\subseteq [\cJ_1]^\al0$ is club.
  $s\in S(C)$ iff for some $N\in C$, $s=N\cap \cI\in S$ and 
  $\eta_s\in \Gen(N)$.
  $S\in S(C^{\cJ_2})$ iff for some 
  $N'\in [\cJ_2]^\al0$ such that $N:=N'\cap J_1\in C$:
  $s=N'\cap \cI\in S$ and 
  $\eta_s\in \Gen(N')$. This is obviously equivalent,
  since $N$ and $N'$ contain the same elements of $I$ and $2^\omega$.
  So $S(C)=S(C^{\cJ_2})$.
  The same argument works with $C\subseteq [\cJ_2]^\al0$
  and $C^{\cJ_1}$. For general $\cJ_1$, $\cJ_2$, apply the argument
  to $\cJ_1$,$\cJ_1\cup \cJ_2$ and $\cJ_2$,$\cJ_1\cup \cJ_2$.
%
 
\item The family $\{\hb(C):\, C \text{ club}\}$ is semi--closed under countable intersections (i.e. if $C_i$ club, $i\in\omega$, then for $C':=\bigcap C_i$ club
  $\hb(C')\subseteq\bigcap \hb(C_i)$). Therefore
  the family $\{B:\, B\supset \hb(C),\, C\text{ club}\}$ is closed
  under countable intersections, and has to contain a minimal element 
  (mod $I$), since $I$ is a ccc--ideal. 
\item 
  Assume, $\hb$ is not truly positive. Wlog $J=H(\chi)$.
   Then for some $C$ club, $B\in I$ Borel: $\hb(C)\subseteq B$.
   Let $C'=\{N\esm H(\chi):\, N\in C, B\in N\}$ club.
   So for any $N\in C'$,  any $Q$--generic over $N$ is not in $B$, so 
   $\hb(C')\subseteq 2^\omega\setminus B$. But 
   $\hb(C')\subseteq \hb(C)\subseteq B$, so $\hb(C')=\emptyset$.
   The rest should be clear.
\end{enumerate}
\end{proof}

\begin{Def}\label{def:prestrue}
A forcing $P$ is called:
\begin{description}  
\item[true positivity preserving] 
   if for all $S,\hb$ truly positive, $\forc_P( \hb\text{ truly pos.})$
\item[true outer measure preserving] 
  if for all $S,\hb$, and $A$ the  true o.m. of $\hb$, 
  $\forc_P( A\text{ is true o.m. of }\hb)$
\end{description}
\end{Def}

These notions do not seem to be equivalent in general (however, they are for
$Q$=Cohen, see Lemma \ref{lem:cohenstrequiv}, and for $Q$=random, provided that
$P$ is weakly homogeneous, see Lemma \ref{lem:weaklyhom}).

Note that true preservation trivially implies properness because of Lemma
\ref{lem:basicstat}(4).

It is clear that true outer measure preserving implies 
true positivity preserving.

\begin{Lem}\label{lem:trueimplpres}
\begin{enumerate}
\item If $P$ is true positivity preserving, then it is positivity preserving.
\item If $P$ is true outer measure preserving, then it is outer measure preserving.
\end{enumerate}
\end{Lem}

\begin{proof}
It is enough to show the following:
For $X$ positive (or: with true outer measure $B$), there
is a $\hb$ truly positive (or: with true outer measure $B$)
s.t. $\{\eta_s:\, s\in S\}\subseteq X$. 
We fix some $\cI\subseteq H(\al1)$ s.t. $\card{\cI}=2^\al0$. 
Wlog $\cJ=H(\al1)$.
\begin{enumerate}
\item 
  For each $N\esm H(\chi)$,
  pick $\eta\in X\cap \Gen(N)$. Clearly, $\hb$ is truly nonempty
  (cf \ref{lem:truebasic}(3)).
\item 
  Let $\beta=2^\al0$.
  As cited in \cite{Kanamori}, $[\cI]^{\al0}$ can be partitioned into $2^\al0$
  many stationary sets, i.e. $[\cI]^{\al0}=\bigcup_{k \in\beta} S_k$.
  Enumerate all positive Borel--subsets of $B$ as $(B_k:\, k\in \beta)$.
  For each $N\esm H(\chi)$, let $k$ be s.t. $N\in S_k$, and 
  pick $\eta\in B_k\cap \Gen^Q(N)$.
  Assume towards a contradiction that the true measure of $\hb$ would 
  be $B'\subset B$, $B_k=B\setminus B'\in \Ip$. If 
  $N\in C\cap S_k$, then $\eta_N\in B_k\cap \hb(C)$, a contradiction.
\end{enumerate}
\end{proof}

As announced, the ``true'' notions are closely related to 
preserving generics:

\begin{Def}\label{def:preservinggenerics}
  $P$ preserves generics, if
  for all $N\esm H(\chi)$, $p\in N$, $\eta\in\Gen^Q(N)$,
  there is a $q\leq p$ $N$--generic s.t. $q\forc \eta\in\Gen^Q(N[G_P])$.
\end{Def}

Notes:\\
1. Instead of for all $N$, we can equivalently say for club many $N$.\\
2. Of course the notion does not depend on $\chi$, provided $\chi$ is 
regular and large enough (in relation to $\card{P}$).\\
3. It is clear that 
preservation of generics
is preserved under composition
(for any Suslin ccc ideal).

Then we get the following:
\begin{Lem}\label{lem:equivtruepres}
  Let $P$ be proper. Then $P$ preserves generics 
  iff it is true positivity preserving.
\end{Lem}


\begin{proof}
  $\rightarrow$:\\
    Assume otherwise, i.e. $\hb$ is truly positive, 
    and $p\forc \hb(\n{C})=\emptyset$. 
    In $V$, $S^*:=\{N\esm H(\chi):\, p,P\in N, N\cap \cI=s\in S, 
    \eta_s\in\Gen(N)\}$ is stationary. (Otherwise,
    the complement of $S^*$ would witness that $\hb$
    is truly empty.) Let $\chi'\gg \chi$,
    $N'\esm H(\chi')$ containing $\hb,S^*,\chi,p,P,
    \n{C}$ such that $N'\cap H(\chi)=N\in S^*$
    (and such that $P$ preserves generics for $N'$,
    if we assume preservation for club many $N$ only).
    So $N'\cap \cI=s\in S$, and
    there is a $q\leq p$ $N'$--generic forcing that 
    $\eta_s\in\Gen^Q(N'[G])$. 
    In $V[G]$, $N'[G]\cap \cI=N'\cap\cI=s$ (since $G$ is $N'$--generic),
    and $N'\cap \n{\cJ}[G]\in \n{C}[G]$  (since $\n{C}\in N'[G]$ club).
    So $\eta_{s}\in \hb(\n{C})$, a contradiction.

  $\leftarrow$:\\
    Assume otherwise, i.e. $N'\esm H(\chi')$, s.t. $p,\eta$ is a counterexample.
    Wlog there is a $\chi\in N$ s.t. $\card{P}\ll \chi \ll \chi'$.
    Let $S:= \{N\esm H(\chi):\, N\text{ is counterexample for }p\text{ and some }\eta\}$
    This set is stationary, since $S\in N'$ and $N'\cap H(\chi)\in S$.

    For each $s\in S$, pick an $\eta_s$ witnessing the counterexample.
    Then $\hb$ is truly positive: If $N\in C\cap S$, then $\eta_N\in \hb(C)$.

    In $V[G]$, let
    $C_\text{gen}:=\{N\esm H^V(\chi):\, G\ N\text{--generic}\}$.
    (Note that the elements of $C_\text{gen}$ are generally not in $V$, only
    subsets of $V$.) 
    $N\esm H^V(\chi)$ just means that $N$ is closed under the Skolem--functions
    of $H^V(\chi)$ (wlog we can also single out a well--order for $H^V(\chi)$, 
    so we just
    need one function), and $G\ N\text{--generic}$ means that 
    for every $D\in N$ such that $D\subseteq P$ is dense,
    $G\cap N\cap D$ is nonempty. Since such $N$ come from simple closure operations,
    $C_\text{gen}$ clearly is club. Therefore also
    $C:=\{\tilde{N}\esm H^{V[G_P]}(\chi):\, G\in \tilde{N}, \tilde{N}\cap V\in C_\text{gen}\}$ is club.
    Therefore, $\hb(C)\neq \emptyset$, i.e. for some $\tilde{N}\esm H(\chi)$, we have:
    $N:=\tilde{N}\cap V\in S\subset V$ and $\eta_N\in\Gen^Q(\tilde{N})$.
    Also, $G$ is $N$--generic, and $N[G]\subseteq \tilde{N}$, so 
    $\eta_N\in\Gen^Q(N[G])$.
 This is a contradiction to the assumption that $\eta_N$ is a counterexample.
\end{proof}

The connection between preservation of true outer measure and 
preservation of generics is a bit more complicated and seems to allow
some variants. Here, we will use the following:

\begin{Def}\label{def:stronglypreserving}
\begin{enumerate}
\item
  $T$ is an interpretation of $\nTs$ w.r.t. $p$, if:
  $T$ is a positive Borel--set, $\nTs$ a $P$--name for a positive Borel--set,
  for all positive Borel--sets $A\subset T$ there is a $p'\leq p$ s.t.
  $p'\forc A\cap \nTs\in\Ip$. 
\item
  $P$ strongly preserves generics, if the following holds:
  For all $N\esm H(\chi)$, $p,T,\nTs\in N$, $T$ an interpretation of $\nTs$
  w.r.t. $p$,
  $\eta\in T\cap \Gen(N)$,
  there is a $q\leq p$ $N$--generic s.t. 
  $q\forc \eta \in \nTs\cap \Gen(N[G_P])$.
\end{enumerate}
\end{Def}

Notes:\\
1. If $p\forc \n{T}''\supset \nTs$, and $T$ is an interpretation
of $\nTs$, then $T$ is an interpretation of $\n{T}''$.\\
2. Again, instead of ``for all $N$'', we can equivalently say 
``for club many $N$'', and the notion does not depend on $\chi$.

\begin{Lem}\label{lem:equivtrueouterpres}
Assume $P$ is proper. Then
\begin{enumerate}
\item 
Preservation of true outer measure implies strong preservation of generics.
\item 
The converse is true provided that 
there are enough interpretations, i.e. the following holds:
If $p\forc \nTs$ is a positive Borel--set, then there 
are $T$, $p'\leq p$ s.t. 
$T$ is an interpretation of $\nTs$ w.r.t. $p'$.
\item 
More generally, we have: $P$ is true outer measure preserving,
if the following holds: If $p\forc \nTs\in\IpB$, then 
there are $T$, $p'\leq p$ such that:\\
(a) $T$ is an interpretation of $\nTs$ wrt $p'$, and\\
(b) for all $N\esm H(\chi)$ s.t. $p,T,\nTs\in N$, for all 
$\eta\in T\cap \Gen(N)$ there is a $q\leq p'$ $N$--generic s.t.
$q\forc \eta\in\nTs\cap \Gen(N[G])$.
\end{enumerate}
\end{Lem}

Note that for $Q$=random (and trivially for $Q$=Cohen),
the additional requirement in (2) is met:
For Cohen, if $p\forc \nTs\in\IpB$, then there are $p'\leq p$ and  a basic open
set  $T$ s.t. $p'\forc T\subseteq \nTs$ (mod $I$).
For random, assume $p\forc \nTs\in \IpB$. 
Wlog $\nTs$ is closed (see note 1 above).
 Let $N\esm H(\chi)$
contain $p,\nTs$, let $G_0\in V$ be $P$--generic over $N$
contain $p$.
Define $T:=\nTs[G_0]$. 
Assume, $A\subseteq T$ is Borel and 
$\Leb(A)>q>0$, $q\in \bQ$.
$\nTs=\bigcap \nTs^n$, where 
$x\in\nTs^n$ iff $\ex y\in\nTs$ s.t. $x\restriction n=y\restriction n$.
The conditions deciding $\nTs$ up to a level 
that is close to the real measure are dense,
i.e. there is an $m\in\omega$, $p'\leq p$ in $N\cap G_0$
such that $p'$ determines $\nTs^m$ 
(i.e. forces it to be $T^m$) and
forces $\Leb(\nTs^m\setminus\nTs)<q/2$.
Let $p''<p'$ be $N$--generic. Then
$p''$ forces that $\Leb(\nTs\Delta T^m)<q/2$,
and $A$ is a subset of $T^m$ of size $>q$,
so $\Leb(A\cap \nTs)>q/2$, i.e. $p$ cannot force
$A\cap \nTs\in I$.

\begin{proof}[Proof of Lemma \ref{lem:equivtrueouterpres}] 
This is similar to the proof of \ref{lem:equivtruepres}.
\begin{enumerate}
\item 
Assume
$p,T,\nTs\in N'\esm H(\chi')$ 
is a counterexample to strong preservation for some $\eta$.
Wlog there is a $\chi\in N$ s.t. $\card{P}\ll \chi \ll \chi'$.
Let $S:= \{N\esm H(\chi):\, N\text{ is counterexample for }p,T,\n{T}
\text{ and some }\eta\}$. Then $S\in N'$, and $N'\cap H(\chi)\in S$,
so $S$ is stationary. For each $N\in S$ let $\eta_N$ be one 
of the counterexamples witnessing that $N\in S$. 

Let $B\subseteq T$
be a true outer measure of $\hb$. So $P$ forces that $B$ is true 
outer measure of $\hb$. 
$B\in \Ip$  (If $C$ is club,
and $N\in C\cap S$, then $\eta_N$ is generic over $N$ since it is 
a counterexample,
so $\eta_N\in \hb(C)$, i.e. $\hb(C)\neq\emptyset$).
So for some $p'\leq p$, $p'\forc B\cap \nTs\in\Ip$. Let $G$ be
$P$--generic, $p'\in G$.

In $V[G]$, define $C_\text{gen}$ and $C$  as in the proof of
Lemma \ref{lem:equivtruepres}. Then $C':=\{\tilde{N}\in C: p'\in \tilde{N}\}$
is club as well. If $\eta_N\in \hb(C')$, then for some 
$\tilde{N}\in C'$, $N=\tilde{N}\cap H^V(\chi)\in S\subset V$,
and $\eta_N\in \Gen(\tilde{N})$. $N[G]\subseteq \tilde{N}$, so
$\eta_N\in \Gen(N[G])$. And since $\eta_N$ is a counterexample,
$\eta_N\notin \nTs[G]$. So $\hb(C')\cap \nTs[G]=\emptyset$,
so the true outer measure of $\hb$ is smaller than $B$, a contradiction.
\item follows from 3.

\item
Assume, $B\supset \hb(C)$ is an outer measure of $\hb$, but in 
$V[G]$, there are $B'$, $C'$ s.t. 
$\hb(C')\subset \hb(C)$, $B'\subset B$, $T':=B\setminus B'\in \Ip$ and 
$B'\supset \hb(C')$. Let this be forced by $p$. So 
according to the assumption we can choose $p'\leq p$, 
$T$ an interpretation of $T'$ w.r.t. $p'$. Wlog $T\subseteq B$.
So $p'$ forces that in
$V[G]$ we get the following picture:\\
\begin{center} 
    \scalebox{0.8}{\begin{picture}(0,0)%
\includegraphics{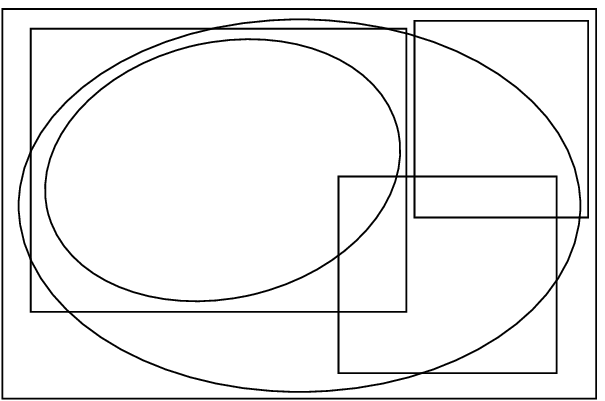}%
\end{picture}%
\setlength{\unitlength}{3947sp}%
\begingroup\makeatletter\ifx\SetFigFont\undefined%
\gdef\SetFigFont#1#2#3#4#5{%
  \reset@font\fontsize{#1}{#2pt}%
  \fontfamily{#3}\fontseries{#4}\fontshape{#5}%
  \selectfont}%
\fi\endgroup%
\begin{picture}(2874,1894)(589,-1493)
\put(676,-1411){\makebox(0,0)[lb]{\smash{\SetFigFont{12}{14.4}{\rmdefault}{\mddefault}{\updefault}{\color[rgb]{0,0,0}$B$}%
}}}
\put(826, 89){\makebox(0,0)[lb]{\smash{\SetFigFont{12}{14.4}{\rmdefault}{\mddefault}{\updefault}{\color[rgb]{0,0,0}$\n{B}'$}%
}}}
\put(3376, 89){\makebox(0,0)[rb]{\smash{\SetFigFont{12}{14.4}{\rmdefault}{\mddefault}{\updefault}{\color[rgb]{0,0,0}$\n{T}'$}%
}}}
\put(3231,-1316){\makebox(0,0)[rb]{\smash{\SetFigFont{12}{14.4}{\rmdefault}{\mddefault}{\updefault}{\color[rgb]{0,0,0}$T$}%
}}}
\put(1001,-811){\makebox(0,0)[lb]{\smash{\SetFigFont{12}{14.4}{\rmdefault}{\mddefault}{\updefault}{\color[rgb]{0,0,0}$\hb(\n{C}')$}%
}}}
\put(1506,-1321){\makebox(0,0)[lb]{\smash{\SetFigFont{12}{14.4}{\rmdefault}{\mddefault}{\updefault}{\color[rgb]{0,0,0}$\hb(C)$}%
}}}
\end{picture}
}
\end{center}
In $V$, $S^*:=\{N\esm H(\chi):\, p',P,T,\nTs\in N, N\cap \cI=s\in S, 
    \eta_s\in\Gen(N)\cap T\}$ is stationary (otherwise, let
$C$ be the complement of $S^*$. Then 
$\hb(C)\cap T=\emptyset$,
so $B$ could not be outer measure of $\hb$).
Let $\chi'\gg \chi$, $N'\esm H(\chi')$ s.t. $S^*, p',P,T,\nTs,\n{C'}\in N'$
and $N:=N'\cap H(\chi)\in S^*$. 
We know that $\eta_N\in T\cap \Gen(N')$
(since $N'\cap 2^\omega=N\cap 2^\omega$),
so by our assumption 
there is a $q\leq p'$ $P$--generic over $N'$ such that
$q\forc \eta_N\in \nTs\cap \Gen(N'[G])$. Also,
$\n{C'}[G]\in N'[G]$, so in $V[G]$, $N'[G]\cap H^{V[G]}(\chi)\in \n{C'}[G]$.
Therefore, $\eta_N\in \hb(\n{C'}[G])$. But $\hb(C')\cap \nTs[G]=\emptyset$,
a contradiction.
\end{enumerate}
\end{proof}

\section{Strong Preservation of Generics for nep Forcings}\label{sec:nep}

In this section, we will prove the following theorem (cf definitions 
\ref{def:preserving} and \ref{def:stronglypreserving}):
\begin{Thm}\label{thm:nep}
If $P$ is nep and Borel outer measure preserving, then 
$P$  strongly preserves generics.
\end{Thm}

\subsection*{About nep Forcings}
Examples for nep (non elementary proper) forcings are Suslin proper forcings
(e.g. Cohen, random, amoeba, Hechler and Mathias)
or Suslin$^+$ forcings (as defined in \cite{tools}, e.g.
Laver, Sacks or Miller).

If you already know what nep forcing is, or you are interested in Suslin proper
forcings only, you can go on directly to the proof of the theorem.  For sake of
completeness, we include a definition of a transitive version of nep here
(which includes e.g.  Laver, Sacks, Miller, see \cite{preservinglittle} for a
proof). In all these cases, in the proof of the theorem $M\vG$ can be
substituted by $M[G]$ (candidates are transitive anyway), and ``ord--collapse''
by ``transitive collapse''.

We assume that the forcing $P$ is defined by formulas $\varphi_{\in
P}(x)$ and $\varphi_\leq(x,y)$, using a real parameter $r_P$.  
Fixing $\ZFCx$, we call $M$ a ``candidate'' if it is a countable transitive
$\ZFCx$ model and $r_P\in M$. So in any candidate, $P^M$ and $\leq^M$ are
defined (but generally not equal to $P\cap M$ or $\leq\cap M$, since 
the definitions do not have to be absolute).

Such a forcing definition $P$ is transitive nep, if 
\begin{enumerate}
\item ``$p\in P$'' and ``$q\leq p$'' are upwards absolute between candidates 
and $V$\\
(i.e. if $M_2\in M_1$, $M_1,M_2$ candidates (or $M_2=V$), and 
$M_1\vD q \in Q$, then $M_2\vD q\in Q$ etc.)
\item In $V$ and all candidates, $P\subseteq H(\al1)$, 
and ``$p\in P$'' and ``$q\leq p$'' are absolute
between the universe and $H(\chi)$ (for large regular $\chi$)
\item For all candidates $M$, $p\in P^M$,
there is a $q\leq p$ s.t. $q\forc (G\cap P^M\text{ is }P^M\text{--gen. over }M)$. (Such a $q$ is called $M$--generic.)
\end{enumerate}



How is this related to proper?
$\ZFCx$ is called normal\label{remark:normal} if for
regular $\chi$ large enough, $H(\chi)\vD\ZFCx$.  We will only be interested
in forcings that are defined with respect to a normal $\ZFCx$.
(Otherwise, if
e.g. $\ZFCx$ contains $0=1$, then every forcing is nep.)
In the normal case, a nep forcing clearly is proper
(consider the transitive collapse of elementary submodels).

In more detail:
Assume $P\subseteq H(\al1)$, $N\esm H(\chi)$ countable,
$i: N\fnto M$ the transitive collapse of $N$. Then $i\restriction P$ 
is the identity, so we have:
$P$ is proper if and only if for all
suitable candidates $M$ and $p\in P^M$ there is a $q\leq p$ $M$--generic,
where suitable means that $M$ is the transitive collapse of an $N\esm H(\chi)$.
Here we allow all candidates, so we get a stronger properness notion.
(Actually, for Theorem \ref{thm:nep} it would be enough to assume 
the properness condition
for internal set forcing extensions of transitive collapses of
elementary submodels only, not for all candidates.)

For Suslin ccc forcings, the
choice of $\ZFCx$ is immaterial, provided that
$\ZFCx$ contains the completeness theorem 
for Keisler--logic. Then any transitive model
of $\ZFCx$ containing the defining real knows 
that $Q$ is a Suslin ccc forcing (see \cite{JdSh:292}).
So we can fix a $\ZFCxQ$ that contains e.g.
the completeness theorem plus the sentences
``there are many regular $\chi$'' and 
``for big regular $\chi$, the completeness theorem
holds in $H(\chi)$''.
It will be implied in the following proof that 
$\ZFCxP$ will include this fixed, finite $\ZFCxQ$.
(And of course we assume that $\ZFCxP$ is normal).

\subsection*{Proof of Theorem \ref{thm:nep}}
The proof is very similar to the proof of ``preserving a little implies
preserving much'' in \cite{Sh:630} (or its version in \cite{preservinglittle}).

From now on, let $M$ be a $P$--candidate, and in $M$: $T$ 
an interpretation of $\nTs$ wrt $p$.

\begin{Def}
    $\es$  is called absolutely $(Q,\h)$--generic
    ($\es\in\aGen(M,p)$), if $\es\in T$ and
    there is a $q\leq p$ $P$--generic over $M$ s.t. (in $V$),
    $p'\forc_P \es\in \nTs\cap\Gen(M\vG)$.
\end{Def}


\begin{Lem}\label{hilfslemma1}
Assume, $P$ is Borel o.m. preserving, 
$M, p, T, \nTs$ as above, 
$M\vD\qb A\in\IpB,\ A\cap T\in \Ip\qe$. 
Then 
$\aGen(M,p)\cap A\in\Ip$.
\end{Lem}

\begin{proof}
Pick (in $M$) a $p'\leq p$ such that $p'\forc A\cap T\cap \nTs\in \Ip$.
Let $q\leq p'$ be $M$--generic. 
In $V[G]$, $\Gen(M[G])$ is $\coI$, and 
$A\cap T\cap \nTs\in \Ip$. Also, $A\cap T$ is o.m. of $(A\cap T)^V$.
So $(A\cap T)^V\cap \nTs\in\Ip$ (otherwise $(A\cap T)\setminus \nTs$
would be the o.m.), so
$X:=\Gen(M[G])\cap V\cap A\cap T\cap \nTs\in \Ip$. 
And clearly $X\subseteq \aGen(M,p)^V\cap A$.
\end{proof}

Assume in $M$,
$2^\card{P}<\chi_1$, $2^{\chi_1}<\chi_2$, $H(\chi_i)\vD\ZFCxP$,
$H(\chi_1)\EQDEF H_1$. Note that for club many $N\esm H(\chi_3))$
($\chi_3$ big enough), the ord--collapse of $N$ is such an $M$.
So it is enough to prove (the obvious analog of) 
strong preservation of generics for these $M$:
If $\es\in \Gen(M)\cap T$, then there is a $q\leq p$ $P$--generic
over $M$ s.t. $q\forc \eta\in\nTs\cap M\vGP$,
i.e. $\es \in \aGen(M,p)$.

Let $R_i$ (in $M$) be the collapse of
$H(\chi_i)$ to $\omega$. 
Let $G_Q\in V$ be a $Q$--generic filter over $M$ s.t.
$\h[G_Q]=\es$, and let
$G_R\in V$ be $R_2$--generic over $M\vGQ$.

\begin{Lem}\label{hilfslemma2} 
$M\vGQ\vGR\vD\qb H_1 \text{ is a (trans.) candidate, }\es\in\aGen(H_1,p)\qe$
\end{Lem}

If this is correct, then Theorem \ref{thm:nep} follows:
Assume, 
$M\vGQ\vGR\vD\qb 
p'\leq  p\ H_1\text{--generic},\ 
p'\forc\es\in\nTs\cap \Gen(H_1[G_P])\qe$.
Let $p''\leq p'$ be $M\vGQ\vGR$--generic. Then $p''$ is $H_1$ generic
and therefore $M$ generic as well
(since $\pow(P)\cap M=\pow(P)\cap H_1$), 
and $p''\forc \es\in\Gen(M\vGP)\cap \nTs$.

\begin{proof}[Proof of Lemma \ref{hilfslemma2}] 
It is clear that 
$H_1$ is a candidate in $M\vGQ\vGR$, and that $\es \in \Gen(H_1)\cap T$.
Assume towards a contradiction, that 
$M\vGQ\vGR\vD\qb \es\notin \aGen(H_1,p)\qe$, Then this is 
forced by some $q\in G_Q$
and $r\in R_2$, but since $R_2$ is homogeneous, wlog $r=1$, i.e.\\
($\ast$) $M\vD\qb q\forc_Q(\h\in T,\, \forc_{R_2} \h\in \Gen(H_1,p) \setminus
\aGen(H_1,p))\qe$.

Now we can construct the following diagram:

\centerline{\xymatrix@!@=0.4cm@M=8pt{
M \ar[rr]^-{R_1} \ar[dr]^-Q_-{G^\otimes_Q} 
&%
&%
{}\save[]+<+2.9cm,0cm>*\txt{$M\la G_{R_1}\ra\EQDEF M_1 \vD \eot\in\aGen(H_1,p)\cap B^M_q$}
\restore%
\ar[rd]^-{R'}_-{\tilde{G}_2}
\\
&
M\la\eot\ra
\ar[ur]^-{R_1/Q}_-{\tilde{G}_1}
&
&%
{}\save[]+<+2.6cm,0cm>*\txt{$M\la\eot\ra\la G_{R_2}\ra=:M_2\vD\eot\notin \aGen(H_1,p)$}
\ar[ll];[]^-{R_2}_-{\tilde{G}_1\ast\tilde{G}_2} 
\restore%
\\
}}

Fix a Borel--set $B^M_q\in M$
s.t. $M\vD\qb\llbracket \h\in B\rrbracket_{\ro(Q)}=q\qe$.
Of course $B^M_q$ is not unique, just unique modulo $I$. 
Such a $B^M_q$ exists, is positive, and we have:\\
$\{\h[G]:\, G\in V\ M\text{--gen},\  q\in G\}
={\ho\omega}\setminus \U\{A^V:\, A\in M, q\forc \h\notin A\}
=\Gen(M)\cap B^M_q$\\
(See e.g. \cite{preservinglittle}). 
$B^M_q\subseteq T$ (mod $I$), 
since $q\forc \h\in T$.
In particular $M\vD\qb B^M_q\cap T\in \IpB\qe$.

Choose $G_{R_1}\in V$ $R_1$--generic over $M$, and let $M_1:=M\langle G_{R_1}\rangle$.
In $M_1$, pick $\eot\in\aGen(H_1,p)\cap B^M_q$ 
(using Lemma \ref{hilfslemma1}), 
so since $\aGen\subseteq \Gen$,
$M_1\vD\qb \ex\, G_Q^\otimes\ Q\text{--gen}/H_1$ s.t. $q\in G_Q^\otimes,
\h[G_Q^\otimes]=\eot \qe$.
This $G_Q^\otimes$ clearly is $M$--generic as well (since 
$M\cap\pow(Q)=H_1\cap\pow(Q)$),
so we can factorize $R_1$ as $R_1=Q\ast {R_1/Q}$ s.t.
$G_{R_1}=G_Q^\otimes\ast \tilde{G}_1$. 

Now we look at the forcing $R_2=R_2^M$ in $M[\eot]$.  
$R_2$ forces that
$R_1$ is countable and therefore equivalent to Cohen forcing.
$R_1/Q$ is a subforcing of $R_1$. Also,
$R_2$ adds a Cohen real. So $R_2$ can be factorized as $R_2=(R_1/Q)\ast R'$,
where $R'=(R_2/(R_1/Q))$. We already have $\tilde{G}_1$  $(R_1/Q)$--generic
over $M[G_Q^\otimes]$, now choose $\tilde{G}_2\in V$ $R'$--generic over $M_1$,
and let $G_{R_2}=\tilde{G}_1\ast\tilde{G}_2$. So $G_{R_2}\in V$ is
$R_2$--generic over $M\la G_Q^\otimes\ra $, $M_2\DEFEQ M\la \eot\ra \la G_{R_2}\ra $.

Let $H_2$ be $H(\chi_2)^{M_1}$. Then 
$H_2\vD\qb p_1\leq p\text{ is }H_1\text{--generic}, p_1\forc \eot\in\Gen(H_1[G_P])\qe$,
and in $M_2$, $H_2$ is a 
candidate. Let in $M_2$, $p_2\leq p_1$ be $H_2$--generic.
Then (in $M_2$), $p_2$ witnesses that $\es\in\aGen(H_1,p)$, 
a contradiction to ($\ast$).
\end{proof}

\section{Preservation for Cohen}\label{sec:cohen}

In this section, let $Q$ be Cohen forcing, i.e. $I$ is the
ideal of meager sets, and $\Gen(N)$ is the set of Cohen reals over $N$.

This is the easiest case: you do not need strong preservation, preservation of
generics itself is iterable; and the proof is a simple modification of the
proof that properness is preserved in a countable support iteration. (This
case could also be seen as a very simple instance of the
general preservation theorem of \cite[XVIII, \S3]{Sh:f},
Case C.) 

We already know that for Cohen, preservation of
Borel positivity is equivalent to preservation of
Borel o.m. The equivalence is also true for the
general preservation notion:

\begin{Lem}\label{lem:cohenstrequiv} 
Preservation of positivity implies preservation of outer measure,
and the same holds for the true version.
\end{Lem}

\begin{proof}
If $A$ is o.m. of $X$, but 
$p\forc( \n B\text{ o.m. of } X,\, A\setminus \n B\in\Ip)$.
Then $A\setminus \n B$ contains a basic open set $D\neq \emptyset$, 
which already exists in $V$. So $p\forc D\cap X\in I$, so
by positivity preservation $D\cap X\in I$, so $A$ cannot be o.m. of $X$.

To show the lemma for the true notion, the same argument works:
Assume, $A$ is true o.m. of $\hb$, and
$p\forc \hb(\n{C}')\cap D\in I$. Then define 
$S^*:=\{s\in S:\, \eta_s\in D\}$, and $\hb^*:=\hb\restriction S^*$.
The usual argument shows that $\hb^*$ is truly positive:
Otherwise, let $C$ be club s.t. $\hb^*(C)=\emptyset$.
Then $C$ witnesses that $A$ is not true o.m. of $\hb$.
On the other hand, $p\forc \hb^*(\n{C}')\in I$,
a contradiction to true positivity preservation.
\end{proof}

\begin{Thm}\label{thm:cohen}
If $(P_i,\n{Q}_i:\, i\in \alpha)$ is a countable support iteration of proper
forcings such that for all $i\in \alpha$, $\forc_{P_i}\n{Q}_i$ preserves
Cohens, then $P_\alpha$ preserves Cohens.
\end{Thm}

\begin{proof}
The successor step is clear, since preservation of generics is 
always preserved by composition.

A real $\eta$ can be interpreted as a function that assigns a natural number to
a sequence of natural numbers. We say $\eta$ is Cohen over a sequence
$(s_0,s_1,\dots,s_{n-1},s_n)$ if $\eta(s_0,\dots,s_{n-1})=s_n$.  
Then $\eta$ is Cohen over $N$ iff for all $f\in N$ there is a $n$ s.t. 
$\eta$ is Cohen over the sequence $f\restriction n$.

Assume, $\alpha=\omega$.
Let $N\esm H(\chi)$ contain $P_\omega$, $p\in P_\omega\cap N$.
Let $\n{f}_i$ and $D_i$ list all $P_\omega$--names for reals and dense sets, resp., 
that are in $N$.

Pick a $p_0\leq p$, $p_0\in N\cap D_0$, s.t. $p_0$ decides $\n{f}_0$
up to a $n_0$ and $\eta$ is Cohen over $\n{f}_0\restriction n_0$.
(This is possible, since inside $N$ we can find an interpretation 
for $\n{f}_0$ and $\eta$ is Cohen over $N$). Then pick
a $q_1\leq p_0\restriction P_1$ $P_1$--generic over $N$ s.t.
$q_1\forc \eta \text{ Cohen over }\Gen(N[G_1])$.

In $V[G_1]$, pick $p_1\leq p_0\in N[G_1]\cap D_1$ s.t. 
$p_1$ proves that $\eta$ is Cohen over $\n{f}_1$ (as above),
and $q_2\leq p_1\restriction P_2$ $P_2$--generic over 
$N[G_1]$ s.t. $q_2\forc \eta \text{ Cohen over }\Gen(N[G_2])$.

Iterating that construction gives us a $q\in P_\omega$
such that $q\restriction P_n\forc q(n)=q_n$, this
$q$ is stronger than $p$ and $N$--generic, and for all
$\n{f}_n$, $q$ forces that $\eta$ is Cohen over $\n{f}_n$.



To prove the theorem for arbitrary $\alpha$, take a sequence $\alpha_i$ $(i\in
\omega)$ cofinal in $\alpha\cap N$. Then do the same as above (however, the
notation and induction gets a bit more complicated, since instead of the $Q_i$
the according quotient forcings have to be used).
\end{proof}

So using the facts that preserving Cohens implies
preserving non--meagerness of arbitrary sets 
(lemma \ref{lem:equivpres2})
and that a nep forcing which 
preserves non--meagerness of Borel--sets preserves Cohens
(Theorem \ref{thm:nep}), we get:

\begin{Cor}
If $(P_i,\n{Q}_i:\, i\in \alpha)$ is a countable support iteration of nep
forcings such that for all $i\in \alpha$, $P_i$ forces that $\n{Q}_i$ preserves
non--meagerness of $V$, then $P_\alpha$ preserves non--meagerness (of all old sets).
\end{Cor}

\section{Preservation for Random}\label{sec:random}

In this section, let $Q$ be random forcing. So $I$ is the
ideal of Lebesgue--Null--sets, and $\Gen(N)$ is the set
of random reals over the model $N$. 

It is clear that preservation of outer measure is equivalent to 
the preservation of the value of $\Leb^*(X)$.

Also the true outer measure is fully 
described by the true outer measure as a real:
Let $\TLeb(\hb):=\min\{\text{Leb}^*(\hb(C)):\, C\text{ club} \}$
(note that  $\TLeb$ really is a minimum).
Then $P$ is true outer measure preserving iff $P$ preserves $\TLeb$.
(This follows from the proof of the next lemma).

\begin{Lem}\label{lem:weaklyhom}
If $P$ is weakly homogeneous
(i.e. if $\varphi$ only contains standard--names,
then $p\forc_P\varphi$ implies $\forc_P\varphi$ ), 
then preserving positivity
implies preserving outer measure, and preserving true
positivity implies preserving true outer measure.
\end{Lem}

\newcommand{\hbs}{{\bar{\eta}^*}}
\newcommand{\nI}{\n{I}}
\begin{proof}
For the ``untrue'' version, this is 
\cite[Lem 6.3.10]{tomekbook}. The same proof works for
true outer measure as well:
Assume that $B$ is a true outer measure of $\hb$,
$\Leb(B)=r_1$
but $p$ forces that $\n{B}'\supseteq \hb(\n{C}')$ and
$\Leb(\n{B}')<r_2<r_1$, $r_2$ rational.
We have to show that there is a truly positive 
$\hbs$ that fails to be truly positive after forcing with $P$.

In $V[G]$ there is a sequence $\nI_n$ of clopen sets 
s.t. $\bigcup \nI_n\supseteq \hb(\n{C}')$ and $\Sigma\Leb(\nI_n)<r_2$.
Let (in $V$) $p_n$, $h(n)$, $I^*_n$ be such that for all $m\leq h(n)$,
$p_n\forc \nI_m=I^*_m \ \&\ \Leb(\bigcup_{m>h(n)}\nI_m)<1/n$.
So $\Leb(\bigcup I^*_m)\leq r_2$. So 
$B\setminus \bigcup I^*_m$ is not null.
Therefore $S^*:=\{s\in S:\, \eta_s\notin \bigcup I^*_m\}$ is
is stationary (otherwise, the complement of $S^*$ 
would witness that $B$ is not the true outer measure).
Define $\hbs:=\hb\restriction S^*$. So $\hbs$ is truly positive.
$p_n\forc \Leb(\bigcup\nI_m\setminus \bigcup I^*_m)
\leq 2/n$, and 
$p_n\forc \hbs(\n{C}')\subseteq \bigcup\nI_m \setminus \bigcup I^*_m$
i.e. $\Leb^*(\hbs(\n{C}'))\leq 2/n$.
So $p_n\forc \TLeb(\hbs)\leq 2/n$. Since this statement does not contain
any names (except standard--names), and $P$ is weakly homogeneous,
$\forc \TLeb(\hbs)\leq 2/n$ for all $n$, i.e. 
$\forc \TLeb(\hbs)=0$. So the truly positive $\hbs$ becomes null
after forcing with $P$.
\end{proof}

For the rest of this section, we will need 
the general iteration theorem of Section 5 of
\cite{tools}, which is cited 
as ``first preservation theorem'' 6.1.B
in \cite{tomekbook}.
It is a simplification of \cite[XVIII,\S3]{Sh:f} Case A.

It uses the following setting: Fix a sequence of increasing arithmetical
two--place relations $R_n$.  Let $R$ be the union of the $R_n$.
Assume $\myC:=\{f:\, f\, R\, g\text{ for some }g\}$ is closed.  $\eta$ covers
$N\esm H(\chi)$, if for every $f\in \myC\cap N$, $f\, R\, \eta$.  We assume
that for every $\eta$ the set $\{f:\, f\, R\, \eta\}$ is closed, and that for
every $N\esm H(\chi)$ there is an $\eta$ covering it.

\begin{Def}\label{def:martinpreserving}
A forcing notion $P$ is tools--preserving, if for all $N\esm H(\chi)$, for all
$p\in P\cap N$, for all $\eta$ covering $N$, for all
$\bar{\n{f}}:=\n{f}_1,\dots,\n{f}_k$ names for elements of $\myC$, and for all
$\bar{f}^*:=f^*_1\dots f^*_k$ interpretations (in $N$) of $\bar{\n{f}}$ under $p$ s.t.
$f^*_i\, R_{n_i}\, \eta$ there is a $q\leq p$ $N$--generic, forcing that $\eta$
covers $N[G]$ and that $\n{f}_i\, R_{n_i}\, \eta$.
\end{Def}

Here, interpretation means that there is an decreasing chain $p=p_0>p_1>\dots$
of conditions s.t. $p_i\forc (\n{f}_1\restriction i=f^*_1\restriction i\, \AND\dots\AND\,
\n{f}_k\restriction i=f^*_k\restriction i)$ (so in particular, $f^*_l\in \myC$).

The general iteration theorem of \cite{tools} says:
\begin{Thm}\label{thm:martin}
Assume, $(P_i,\n{Q}_i:\, i<\alpha)$ is a countable support iteration
of proper, tools--preserving forcings. Then $P_\alpha$ is 
tools--preserving.
\end{Thm}

In the case of random, we list the clopen subsets of $2^\omega$ 
as $(I_i:\, i\in\omega)$, and interpret a function $f$ as 
a sequence of clopen sets. We let $\myC:=\{f:\, \forall i\, \text{Leb}(I_{f(i)})<2^{-i}\}$,
and define $f\, R_n\, \eta$ by: for all $l>n$, $\eta\notin I_{f(l)}$.
Then $\eta$ covers $N$ iff $\eta$ is random over $N$ (see eg.
\cite{tomekbook} or \cite{tools}).

\begin{Lem}\label{lem:trivial}
For random, the following are equivalent:
\begin{enumerate}
\item
$P$ is tools--preserving.
\item 
$P$ is tools--preserving for $k=1$ and $n_1=0$.
\item
$P$ is strongly preserving randoms.
\end{enumerate}
\end{Lem}

\begin{proof}
$2\rightarrow 1$:
So we have given
$f^*_1\dots f^*_k$ an interpretation of $\n{f}_1\dots,\n{f}_k$,
$p$, $N$, $\eta$.
Let $n^*:=\max(k,n_1,\dots,n_k)$.
Define $g^*$ s.t. $I_{g^*(m)}=\bigcup_{i=1\dots k} I_{f^*_i(n^*+m)}$,
and $\n{g}$ s.t. $I_{\n{g}(m)}=\bigcup_{i=1\dots k} I_{\n{f}_i(n^*+m)}$.
So for all $m$, $\Leb(I_{g^*(m)})<k 2^{-(n^*+m)}$, so $p\forc \n{G}\in\myC$,
and $g^*$ is an interpretation of $\n{g}$.
Also, for all $m$, $\eta\notin I_{g^*}$, i.e. $\eta R_0 g^*$.
Let $p'\leq p$ s.t. $p'\forc f^*_i\restriction n^*=\n{f}_i\restriction n^*$.
So by 2, there is a $q\leq p'$ $N$--generic s.t. $q$ forces that
$\eta$ is random over $N[G]$ and that $\eta R_0 \n{g}$.
So for all $m>n^*$, $\eta\notin I_{\n{f}_i(m)}$.
And for $n_i\leq m < n^*$, $I_{\n{f}_i(m)}=I_{f^*_i(m)}\notni \eta$,
so $q$ forces that $\eta R_{n_i} f_i$.

$2 \rightarrow 3$:
It is enough to show that the assumptions for
\ref{lem:equivtrueouterpres}(3) are met.
So let $p\forc\nTs\in \IpB$. 
Then (in $V[G]$) there is a closed subset $\n{A}$ of $\nTs$ that is positive.
Let $\n{X}$ be the family of all clopen supersets of $\n{A}$.
Clearly $\{\Leb(I):\, I\in \n{X}\}$ is dense in the interval $[\Leb(\n{A}),1]$.
So we can find a decreasing sequence $\n{I}^n$ of clopen supersets of $A$ s.t.
$2^{-n}<\Leb(\n{I}^n\setminus \n{A})<2^{-(n-1)}$.
Let $\n{f}$ code the sequence $\bar{\n{I}}^n:=\n{I}^n\setminus \n{I}^{n+1}$.
Then $\n{f}\in \myC$. 
Now in $V$, pick any $N'\esm H(\chi')$ containing $p,\n{f}$
and let $G\in V$ be $N'$--generic. Then $f^*=\n{f}[G]$
is an interpretation of $\n{f}$ (in the sense
of tools--preservation). 
Let $p'\leq p$ force this, and force a value to $\n{I}^0$.
$\n{f}$ defines a sequence of clopen sets $\bar{I}^n$.
Let $T:=I_0\setminus \bigcup \bar{I}^n$. Then $\Leb(T)>0$,
and $T$ is an interpretation of $\nTs$.
Let $N\esm H(\chi)$ contain $T,\nTs,\dots$ and let
$\eta\in T\cap \Gen(N)$.
Then $\eta R_0 f$, so there is a $q\leq p$ 
$N$--generic
forcing that $\eta\in\Gen(N[G])$ and that $\eta R_0 \n{f}$.
Since $\eta \in I^0$, $q$ forces that $\eta\in\nTs$.

$3 \rightarrow 2$:
Given $f^*$ and $f$, define $T:=\cap 2^\omega\setminus I_{f^*(m)}$,
and the same for $\nTs$ and $\n{f}$.
Then $T$ is an interpretation of $\nTs$.

\end{proof}

Using Theorem \ref{thm:nep}, the fact that strong preservation implies
preservation (see e.g. Lemma \ref{lem:equivpres2})
and the last lemma we get:

\begin{Cor}
Assume, $(P_i,\n{Q}_i:\, i<\alpha)$ is a countable support iteration
of nep forcings s.t. for all $i$, 
$P_i$ forces that $\n{Q}_i$ preserves Lebesgue--positivity of $V$.
 Then $P_\alpha$ preserves Lebesgue--positivity (of all old sets).
\end{Cor}
\newpage
\section*{The Diagram of Implications}\label{sec:diagramimpl}
\addcontentsline{toc}{section}{The diagram of implications}

For general Suslin ccc ideals we get:\\[3ex]
\label{diagram:implications}\centerline{\xymatrix@M+=2ex{
*+\txt{pres. true\\outer measure}
\ar[ddd]^{\text{\ref{lem:trueimplpres}}}
\ar[rr]
&
&
*+\txt{pres. true\\positivity}\ar@{<->}[d]^{\text{\ref{lem:equivtruepres}}}
\\
&
*+\txt{strongly\\pres. generics}\ar@<-0.5ex>@{-->}[ul]_{\text{many interpret.: \ref{lem:equivtrueouterpres}}}
\ar@<-0.5ex>[ul];[]
\ar[r]
&
*+\txt{preserving\\generics} \ar[d]
\\
&
&
*+\txt{pres. many\\generics} \ar@{<->}[d]^{\text{\ref{lem:equivpres2}}} 
\\ 
*+\txt{preserving\\ outer measure} \ar[d]
\ar[rr]
& 
&
*+\txt{preserving\\positivity} \ar[d]
\\
*+\txt{pres. Borel\\outer measure} \ar@{<->}[d]
\ar[rr]
\ar@/_1pc/@{-->}[uuur]|(0.65){P \text{ nep: \ref{thm:nep}}} 
& 
&
*+\txt{pres. Borel\\positivity} 
\ar@<-0.5ex>[d]\ar@{-->}@<-0.5ex>[d];[]_{P \text{ Borel--hom: \ref{lem:borelhom}}}
\ar@/_3pc/@{-->}[uuu]|(0.55){P \text{ nep: \ref{thm:nep}}} 
\\
*+\txt{pres. outer\\measure of $V$}
& 
&
*+\txt{preserving\\positivity of $V$}
\ar@{-->}[ll]_{\text{o.m. of }V\ \emptyset\text{ or }2^\omega}
\\
}}
\vspace{2ex}
Preservation of (Borel) positivity and outer measure is defined in
\ref{def:preserving}, the true notions in \ref{def:prestrue},
and (strong) preservation of generics in 
\ref{def:preservinggenerics} and \ref{def:stronglypreserving}.

For ``many interpretations'' and ``$P$ nep'', see 
\ref{lem:equivtrueouterpres}(2) and Section \ref{sec:nep}, resp.
\newpage
In the special case of random and Cohen we get:\\[3ex]
\begin{minipage}[t]{0.59\textwidth}
$Q$=random (i.e. $I$=Lebesgue--null):\\
\xymatrix@M+=2ex{
*+\txt{tools--pres.\\(iterable)$\IFF$\\ true o.m. pres$\IFF$\\str. pres. gen.}
\ar[d]
\ar@<-0.5ex>[r] \ar@{-->}@<-0.5ex>[r];[]_{P \text{ hom: \ref{lem:weaklyhom}}} 
&
*+\txt{ true pos. pres}\ar[d]
\\
*+\txt{outer measure\\preserving}
\ar@<-0.5ex>[r] \ar@{-->}@<-0.5ex>[r];[]_{P \text{ hom: \ref{lem:weaklyhom}}} 
&
*+\txt{pres. pos.}
\\
&
{}\save[]+<0cm,-2ex>*\txt{pres. pos. of $V$ $\IFF$\\pres. Borel o.m.}
\ar[u];[]
\ar@/^1pc/@{-->}[luu]|(0.3){\txt{\footnotesize $P$ nep: \ref{thm:nep}}} 
\restore \\
}
\\[2ex]
For the definition of tools--preserving,
see \ref{def:martinpreserving}.
``$P$ hom'' means 
$P$ is weakly homogeneous, see \ref{lem:weaklyhom}.
\end{minipage}\hfill
\begin{minipage}[t]{0.40\textwidth}
$Q$=Cohen (i.e. $I$=meager):\\
\xymatrix@M+=2ex{
*+\txt{true o.m. pres.$\IFF$\\pres. gen. (iterable)}\ar[d]
\\
*+\txt{o. m. pres.\IFF\\ pos. pres.}\ar[d]
\\
*+\txt{pres. pos. of $V$ $\IFF$\\pres. Borel o.m.}
\ar@/^6pc/@{-->}[uu]|(0.3){\txt{\footnotesize $P$ nep: \ref{thm:nep}}} \\
}
\end{minipage}%

\bibliographystyle{amsalpha}
\bibliography{logik,listb}

\end{document}